\documentstyle[11pt]{article}
\begin{document}

\title{Systems of conservation laws of Temple class,
        equations of associativity and linear congruences in $P^4$}

\author{{\Large Agafonov S.I. } \& {\Large Ferapontov E.V. }\\
\\
    Department of Mathematical Sciences \\
    Loughborough University \\
    Loughborough, Leicestershire LE11 3TU \\
    United Kingdom \\
    e-mails: \\
    {\tt S.Agafonov@lboro.ac.uk} \\
    {\tt E.V.Ferapontov@lboro.ac.uk}
}
\date{}
\maketitle

\newtheorem{theorem}{Theorem}
\newtheorem{proposition}{Proposition}
\newtheorem{lemma}{Lemma}

\pagestyle{plain}

\maketitle

\begin{abstract}
We propose a geometric correspondence between (a) linearly degenerate
systems of conservation laws with
rectilinear rarefaction curves and (b) congruences of lines in projective
space whose developable surfaces are planar
pencils of lines. We prove that in $P^4$ such congruences are necessarily
linear. Based on the
results of Castelnuovo,
the classification of three-component systems is obtained, revealing a close
relationship of
the problem with  projective geometry of the Verones\'e variety $V^2 \subset
P^5$ and 
the theory of associativity equations of  two-dimensional topological field
theory.

\bigskip

Subj. Class.: Differential Geometry, Partial Differential Equations.



\end{abstract}

\section{Introduction}

Hyperbolic systems of conservation laws
\begin{equation}
u^i_t=f^i(u)_x=v^i_j(u)  u^j_x, ~~~
v^i_j=\frac{\partial f^i}{\partial u^j},~~~i=1,...,n,
\label{cons}
\end{equation}
naturally arise in a variety of physical applications and are known to
possess a rich mathematical and geometric structure
\cite{Jef}, \cite{Lax}, \cite{Sevennec}, \cite{DN}, \cite{Tsarev}.
It was observed recently that many constructions of the theory of systems of
conservation laws 
are, in a sense, parallel to that of the projective theory of congruences.
The correspondence proposed in  \cite{Fer1} and \cite{Fer2} associates with
any system (\ref{cons})
an $n$-parameter family of lines
\begin{equation}
y^i=u^i\ y^0-f^i(u), ~~~ i=1,..., n
\label{I2}
\end{equation}
in $(n+1)$-dimensional projective space $P^{n+1}$ with affine coordinates
$y^0,..., y^n$.
In the case $n=2$ we obtain a two-parameter family, or a congruence of lines
in $P^3$. In the 19th century
the theory of congruences was one of the most popular chapters of classical
differential geometry
(see, e.g., \cite{Finikov}).
We keep the name``congruence" for any $n$-parameter family of lines
(\ref{I2}) in $P^{n+1}$.

It turns out that the basic concepts of the theory of systems
of conservation laws, such as  shock and rarefaction curves, Riemann
invariants, 
reciprocal transformations, linearly degenerate systems and systems of
Temple class \cite{Temple}
acquire a clear and simple projective interpretation when reformulated in
the 
language of the theory of congruences. For instance, this correspondence
enabled the classification of systems of Temple class to be reduced to a
much
simpler geometric problem of the classification of congruences with either
planar or 
conical developable surfaces. In particular, some of the results of
\cite{Temple}
became intuitive geometric statements about families of lines in projective
space.
Another application of the correspondence proposed was the construction of
the
Laplace and L\'evy transformations of hydrodynamic type systems in Riemann
invariants
\cite{Fer3}, \cite{Fer5} which, on the geometric level, have been a subject
of 
extensive research in  projective differential geometry.

\noindent {\bf Remark.} It should be emphasized that the correspondence
between  systems (\ref{cons}) and
congruences in $P^{n+1}$ is not one-to-one: some ``degenerate" congruences
are to be excluded. Indeed, let
\begin{equation}
y^i=g^i(u) y^0-f^i(u),~~~u=(u^1,...,u^n)
\label{congr}
\end{equation}
be an arbitrary $n$-parameter family of lines in $P^{n+1}$; notice that
$g^i(u)$, as well as $f^i(u)$, may not
happen to be fuctionally independent.
Associated with such a congruence is a system of conservation laws
\begin{equation}
g^i(u)_t=f^i(u)_x, ~~~ i=1,..., n,
\label{coneq}
\end{equation}
which, for functionally dependent $g^i(u)$, is not in  Cauchy normal form.
System (\ref{coneq}) can be transformed to the Cauchy normal form provided
the characteristic polynomial
\begin{equation}
{\rm det}\left(\lambda \frac{\partial g^i}{\partial u^j} - \frac{\partial
f^i}{\partial u^j}\right)
\label{nondegenerate}
\end{equation}
is not identically zero which, on the geometric level, is equivalent to the
requirement 
that the lines (\ref{congr})
do not belong to a  hypersurface in $P^{n+1}$. Hypersurfaces in $P^{n+1}$
carrying 
$n$-parameter families of lines are  interesting
 in their own. For $n=2$ these are planes. In the case $n=3$ these are
either one-parameter families of planes
or three-dimensional quadrics
\cite{Segre1}. For $n=4$ among obvious examples are two-parameter families
of planes or one-parameter families of
three-dimensional quadrics -- see \cite{Togliatti}, \cite{Rogora} for the
classification results.
In what follows we consider nondegenerate {\it hyperbolic} congruences only,
which means that the characteristic polynomial
(\ref{nondegenerate})
is not zero identically  and its roots are real and pairwise distinct. Any
such congruence can be parametrized in the form
(\ref{I2}). 

\medskip

Let $\lambda^i(u)$ be the eigenvalues of the matrix $v^i_j$ of system
(\ref{cons}), assumed real and pairwise distinct.
Let $\xi^i(u)$ be the corresponding eigenvectors: $ v \xi^i=\lambda^i
\xi^i$. In this paper we investigate and classify systems of
conservation laws which simultaneously satisfy the following two properties:

(a) The integral trajectories of  the eigenvectors $\xi^i$ (called the
rarefaction curves of system (\ref{cons})) are
straight lines in coordinates $u^1, ..., u^n$. This condition was introduced
by Temple in \cite{Temple}.

(b) The eigenvalues $\lambda^i$ are constant along rarefaction curves of the
i-th family. Such systems are known as linearly degenerate.

Systems (\ref{cons}) satisfying both these conditions will be called {\bf
T-systems} for short.
In section 2 we include the necessary information about systems of
conservation laws and recall the main results of
Temple \cite{Temple} and our recent work \cite{Fer1}, \cite{Fer2}. These
results imply the correspondence between
$n$-component T-systems and congruences
in $P^{n+1}$ whose developable surfaces are planar pencils of lines. One can
readily establish that for
$n=2$ such congruences are linear (that is, defined by two linear equations
in Pl\"ucker coordinates)
and
consist of all lines intersecting two fixed skew lines in $P^3$. Since any
two linear congruences in $P^3$ are projectively equivalent,
there exists essentially a unique two-component T-system. Recall that
projective 
transformations of congruences
(\ref{I2}) correspond to ``reciprocal transformations'' of systems
(\ref{cons}) which preserve the T-property:
see section 3 for the details.

{\bf Example 1.} Consider the wave equation
\begin{equation}
f_{tt}-f_{xx}=0.
\label{wave}
\end{equation}
Introducing the variables $a=f_{xx}, \ b=f_{xt}$, we readily rewrite
(\ref{wave}) 
as a linear two-component system of conservation laws
\begin{equation}
a_t=b_x, ~~~~~ b_t=a_x
\label{cons1}
\end{equation}
which is obviously a T-system (any linear system of conservation laws is a
T-system since its eigenvalues and
eigenvectors are constant). The corresponding congruence (\ref{I2})
\begin{equation}
y^1=ay^0-b, ~~~~~ y^2=by^0-a
\label{I22}
\end{equation}
consists of all lines intersecting the two skew lines  $y^0=1, \ y^1=-y^2$
and $y^0=-1, \ y^1=y^2$.

{\bf Example 2.} Consider the Monge-Amp\`ere equation
\begin{equation}
f_{xt}^2-f_{xx}f_{tt}=1.
\label{Monge}
\end{equation}
Introducing the variables $a=f_{xx}, \ b=f_{xt}$ \cite{Mokhov2}, we readily
rewrite (\ref{Monge})
as a two-component system of conservation laws
\begin{equation}
a_t=b_x, ~~~~~ b_t=\left(\frac{b^2-1}{a}\right)_x
\label{cons2}
\end{equation}
which proves to be a T-system. The corresponding congruence
\begin{equation}
y^1=ay^0-b, ~~~~~ y^2=by^0-\frac{b^2-1}{a}
\label{I23}
\end{equation}
consists of all lines intersecting the two skew lines  $y^1=1, \ y^0=y^2$
and $y^1=-1, \ y^0=-y^2$.

Since congruences (\ref{I22}) and (\ref{I23}) are projectively equivalent,
the corresponding 
systems (\ref{cons1}) and
(\ref{cons2}) are reciprocally related, thus providing a linearization of
the nonlinear Monge-Amp\`ere equation
(\ref{Monge}) (which, of course, is not a new result).

\bigskip

Our main result is the classification of three-component T-systems or, in
geometric language, congruences in $P^4$
whose developable surfaces are planar pencils of lines. The main example
which motivated our research comes from the theory
of equations of associativity  of two-dimensional topological field theory.

{\bf Example 3.} Let us consider the Monge-Amp\`ere type equation
\begin{equation}
f_{ttt}=f_{xxt}^2-f_{xxx}f_{xtt},
\label{ass}
\end{equation}
known as the WDVV or the associativity equation, which was thoroughly
investigated by Dubrovin in \cite{Dub}.
Introducing the variables $a=f_{xxx}, \ b=f_{xxt}, \ c=f_{xtt}$
\cite{Mokhov1},  we readily rewrite (\ref{ass})
as a three-component system of conservation laws
\begin{equation}
a_t=b_x, ~~~~~ b_t=c_x, ~~~~~ c_t=(b^2-ac)_x
\label{cons3}
\end{equation}
which was observed to be a T-system in \cite{Fer1}. The corresponding
congruence  in $P^4$
\begin{equation}
y^1=ay^0-b, ~~~~~ y^2=by^0-c, ~~~~~ y^3=cy^0-b^2+ac
\label{I24}
\end{equation}
 coincides with the set of trisecant lines of the Veronese variety
 projected from $P^5$ into $P^4$ (see section 3).
In this sense the projected Veronese variety plays the role of the focal
variety of the congruence (\ref{I24}).
As follows from the classification result presented below, this example is,
in a sense, generic.

We prove that  congruences in $P^4$ whose developable surfaces are planar
pencils of lines are necessarily linear
(that is, defined by three linear equations in the Pl\"ucker coordinates).
In the parametrisation (\ref{I2}) the Pl\"ucker coordinates of a congruence
in $P^4$
are 
$$
1, \ ~~ u^1, \ u^2, \ u^3, \ ~~~ f^1, \ f^2, \ f^3, \ ~~~ u^1f^2-u^2f^1, \
u^1f^3-u^3f^1, \ u^2f^3-u^3f^2.
$$
Linear congruences are characterized by three linear relations among them
$$
\alpha +\alpha _iu^i+\beta _if^i+\alpha _{ij}(u^if^j-u^jf^i)=0,
$$
where $\alpha, \ \alpha_i, \ \beta_i, \ \alpha_{ij}$ are arbitrary
constants. Solving these equations
for $f^1, \ f^2, \ f^3$, we arrive at the general formula for the fluxes of
three-component T-systems.
Notice that the congruence (\ref{I24}) is linear.

{\bf Remark.} Given a congruence in $P^4$ whose developable surfaces are
planar pencils of lines, its Pl\"ucker image
is a three-dimensional submanifold of the Grassmanian $G(1, 4)$ covered by a
two-parameter family of lines,
the images of planar pencils. Moreover, there are three lines passing
through each point of this submanifold.
Algebraic threefolds covered by lines were recently classified by Mezzetti
and Portelli in \cite{Mezzetti1}. It was demonstrated
that threefolds with three lines through each point are intersections of
$G(1, 4)$ with a $P^6$, that is,
Pl\"ucker images of linear congruences. In this paper we approach the
classification problem from the point of view of
local differential geometry,
without imposing any additional algebro-geometric restrictions.
It turns out, however, that our local differential-geometric assumptions
(namely, that developable surfaces are planar pencils)
already imply the algebraizability.

Unlike the case of $P^3$, the proof of the linearity of these congruences in
$P^4$ requires a long computation
bringing a certain exterior differential system into involutive form, which
is carried over to the Appendix (notice that the linearity does not
necessarily hold in $P^5$ as
simple examples from section 5
show). Once the linearity is established, one can make use of the results of
Castelnuovo 
\cite{Castelnuovo} who classified  linear congruences in $P^4$. He found six
projectively different types thereof,
thus providing a list of six three-component T-systems which are not
reciprocally related.
 Below we list them  as  scalar third-order
Monge-Amp\`ere type equations. They assume the form (\ref{cons}) in the
variables $a=f_{xxx}, \ b=f_{xxt}, \ c=f_{xtt}$.
As systems of conservation laws, they differ by a number of Riemann
invariants they possess (see section 2 for the definitions).
Geometrically, the existence of a Riemann invariant implies the reducibility
of the focal variety of
the corresponding congruence: if a T-system possesses $k$ Riemann
invariants, 
the focal variety contains $k$ linear subspaces of codimension two.

\begin{theorem}
\label{MainList}
 Any strictly hyperbolic three-component T-system can be reduced by a
reciprocal transformation to one
from the following list.

\medskip

\noindent
{\bf I.} T-systems which possess no Riemann invariants,
\begin{equation}
f_{xxx}f_{ttt}-f_{xxt}f_{ttx}=1
\label{eq1}
\end{equation}
and
\begin{equation}
f_{xxt}^2+f_{xtt}^2-f_{xxx}f_{xtt}-f_{ttt}f_{xxt}=1.
\label{eq2}
\end{equation}
The focal varieties of the corresponding congruences  are non-singular
projections of the Verones\'e variety into $P^4$.
The congruences consist of the trisecant lines of these projections. Notice
that  there are two  different
projections which are not equivalent over the reals.

\medskip

\noindent {\bf II} T-systems which possess  one Riemann invariant,
\begin{equation}
f_{xxx}f_{ttt}-f_{xxt}f_{ttx}=0
\label{eq3}
\end{equation}
and
\begin{equation}
f_{xxt}^2+f_{xtt}^2-f_{xxx}f_{xtt}-f_{ttt}f_{xxt}=0.
\label{eq4}
\end{equation}
The corresponding  
focal varieties are reducible and consist of a cubic scroll and a plane
which intersects the cubic scroll along its directrix.
Notice that equations (\ref{eq1}) and (\ref{eq3}) are related to (\ref{eq2})
and (\ref{eq4}) by a complex change of variables
$x\to (x+t)/\sqrt 2, \ t\to i(x-t)/\sqrt 2$.

\medskip

\noindent {\bf III} T-system with two Riemann invariants,
\begin{equation}
f_{xtt}^2-f_{xxt}f_{ttt}=1
\label{eq5}
\end{equation}
which  reduces to the Monge-Amp\`ere equation (\ref{Monge}) for  $\tilde f
=f_t$.
The corresponding focal variety consists of a two-dimensional quadric and
two planes
which intersect the quadric along rectilinear generators of different
families.

\medskip

\noindent {\bf IV} T-system with three Riemann invariants,
\begin{equation}
f_{ttt}-f_{xxt}=0.
\label{eq6}
\end{equation}
The corresponding focal variety consists of three planes.
\end{theorem}

\noindent We discuss the geometry of these examples in some more detail in
section 3.

\noindent {\bf Remark.} Equation (\ref{eq1})  was discussed by Dubrovin in
\cite{Dub}. 
As shown in \cite{Fer6}, after the  transformation
$
\tilde x=t, ~ \tilde t= f_{xx}, ~ \tilde f_{\tilde x \tilde x}=-f_{xt}, ~
\tilde f_{\tilde x \tilde t}=x, ~
\tilde f_{\tilde t \tilde t}=f_{tt},
$
it takes the form (\ref{ass}):
$
\tilde f_{\tilde t \tilde t \tilde t}=\tilde f_{\tilde x \tilde x \tilde
t}^2-
\tilde f_{\tilde x \tilde x \tilde x}\tilde f_{\tilde x \tilde t \tilde t}.
$
Notice that this is not a contact transformation. Geometrically, equations
(\ref{ass}) and (\ref{eq1}) correspond to projectively equivalent
congruences (see section 3).
Equation (\ref{eq3}) was discussed in \cite{Fer7} and \cite{Strachan}. The
classification of third order equations
of Monge-Amper\'e type was given in \cite{A}.

We would like to conclude this introduction by formulating two conjectures
about the structure of
congruences in $P^{n+1}$
whose developable surfaces are planar pencils of lines.

\medskip

\noindent 1. The focal varieties of such congruences are algebraic
(possibly, reducible and singular).

\medskip 

\noindent 2. The intersection of the focal variety with a developable
surface (which is a planar pencil of lines)
consists of a point (the vertex of the pencil) and a plane curve of degree
$n-1$.

For $n=2$ this is obvious. For $n=3$ it follows from the results presented
above. Both conjectures are true for general linear
congruences in $P^{n+1}$ (see section 4) and congruences arising from the
completely exceptional Monge-Amp\`ere type equations
(see section 5). As readily follows from the discussion in section 2, the
focal varieties have
codimension  two and contain n-parameter families of plane curves (which are
conics for $n=3$). This shows that the problem in question is actually
algebro-geometric.

\section{T-systems from the point of view of the projective theory of
congruences.} 

In this section we give a brief review of the necessary material from
\cite{Fer1}, \cite{Fer2} and \cite{Temple}.
Let $\lambda^i(u)$ be the eigenvalues of the matrix $v^i_j(u)$
(called the characteristic velocities of system
(\ref{cons})), assumed  real
and pairwise distinct. Let ${\bf \xi }_i=(\xi ^1_i,...,\xi^n_i)^T$ be the
corresponding eigenvectors:
$v{\bf \xi }_i=\lambda^i{\bf \xi }_i$, or, in components, $v^s_k \xi
^k_i=\lambda ^i \xi ^s_i.$ We denote by
$L_i=\xi^k_i\frac{\partial}{\partial u^k}$ the Lie derivative in the
direction of  ${\bf \xi}_i$.
It is convenient to introduce the expansions
\begin{equation}
[L_i,L_j]=c^k_{ij}L_k,~~~c^k_{ij}=-c^k_{ji}.
\label{comm}
\end{equation}
In the theory of  hydrodynamic type systems  {\it rarefaction curves} play
crucial role. Recall that rarefaction curves
are integral curves of the eigenvectors ${\bf \xi }_i$.
Thus, there are $n$ families of rarefaction curves, and for any point in
$u$-space there is exactly one rarefaction curve
from each family passing through it.
Due to
the correspondence (\ref{I2}), a curve in $u-$space
defines a ruled surface, i.e., a one-parameter family of lines in $P^{n+1}.$

\begin{theorem}\cite{Fer1}
Ruled surfaces defined by rarefaction curves of the i-th family are
developable, i.e., 
their rectilinear generators are tangential to a curve.
This curve can be parametrized in the form
\begin{equation}
y^0=\lambda ^i,~~y^1=u^1\lambda ^i -f^1(u),~ ..., ~ y^n=u^n\lambda ^i
-f^n(u),
\label{cuspidal}
\end{equation}
where $u$ varies along the rarefaction curve.
\label{rare}
\end{theorem}
The curve (\ref{cuspidal}) constitutes the singular locus of the developable
surface which is called its {\it cuspidal edge}.
The collection of all cuspidal edges corresponding to rarefaction curves of
the i-th family defines the
 so-called focal hypersurface $M_i \subset P^{n+1}$. In our case parametric
equations of $M_i$ coincide with (\ref{cuspidal}),
 where $u$
is now allowed to take all possible values. By a construction, each line of
the 
congruence (\ref{I2}) is tangential to the focal
hypersurface $M_i $. The idea of focal hypersurfaces is obviously borrowed
from optics: thinking of the lines of a
congruence as the rays of light, one can intuitively imagine focal
hypersurfaces as the locus in $P^{n+1}$ where the light
concentrates (this explaines why in German literature focal hypersurfaces
are called 'Brennfl\"achen', which can be translated
as 'burning surfaces').

Since the system of conservation laws (\ref{cons}) is strictly hyperbolic,
there are precisely $n$ developable surfaces passing through a
line of the congruence (\ref{I2}), and each line is tangential to $n$ focal
hypersurfaces. 

In the theory of weak solutions of  systems (\ref{cons})  {\it shock curves}
play fundamental role.
The shock curve with the vertex in $u_0$ is the set of points in the
$u$-space such that
\begin{equation}
\sigma(u^i-u^i_0)+f^i(u)-f^i(u_0)=0,~~i-1,...,n,
\label{shock}
\end{equation}
for some function $\sigma(u,u_0).$
For any $u$ on the shock curve the discontinuous function
$$
\begin{array}{l}
u(x,t)=u_0, ~~x\le \sigma t,\\
u(x,t)=u, ~~x\ge \sigma t,
\end{array}
$$
is a weak solution of (\ref{cons}).
Shock curves, like  rarefation curves, define special ruled surfaces of the
congruence (\ref{I2}), the geometry of which was
clarified in \cite{Fer1}.

Lax showed that a shock curve with the vertex in a generic point $u_0$
splits into $n$ branches, the $i$th branch being
 $C^2-$tangent of the associated rarefaction curve of the $i$th family
passing through $u_0.$
 
 As pointed out by a number of authors, there are situations when shock
curves coincide with their associated
 rarefaction curves. Systems with coinciding shock and rarefaction curves
were studied by Temple \cite{Temple}.
 His main result can be formulated as follows.
 \begin{theorem} \cite{Temple}
Rarefaction curves of the $i$th family coincide with the associated
branches of the shock curve if and only if either
 
 1) every rarefaction curve of the $i$th family is a straight line in the
$u-$space 
 
 or 
 
 2) the characteristic velocity $\lambda ^i$ is constant along rarefaction
curves of the $i$th family,
 $$
 L_i(\lambda ^i)=0.
 $$
 \label{Te}
 \end{theorem}
Systems satisfying the condition 2 are known as {\it linearly degenerate}.
Both these condition
have a very natural geometric interpretation.
\begin{theorem}\cite{Fer1}
Rarefaction curves of the $i$th family are straight lines if and only if the
associated developable surfaces are planar,
that is, their cuspidal edges are plane curves.
\label{line}
\end{theorem}
\begin{theorem}\cite{Fer1}
The characteristic velocity $\lambda ^i$ is linearly degenerate if and only
if the associated developable surfaces  are conical,
that is, their generators meet in a point. The corresponding focal
hypersurface $M_i$ degenerates into a submanifold of codimension two.
\label{conic}
\end{theorem}
As demonstrated in \cite{Fer1} and \cite{Fer2}, theorems \ref{line} and
\ref{conic} provide  an elementary geometric proof
of Theorem \ref{Te}.

In what follows, we consider systems (\ref{cons}) which simultaneously
satisfy both conditions of Theorem \ref{Te}, namely,
all rarefaction curves are rectilinear (in u-space), and all eigenvalues
$\lambda ^i$ are linearly degenerate. Systems of this
type naturally arise in the theory of equations of associativity of 2D
topological field theory
\cite{Dub}  (see examples below). We will
call them {\bf T-systems} for short. In  view of Theorems \ref{line} and
\ref{conic}, developable surfaces of the corresponding
congruence (\ref{I2}) are planar and conical simultaneously, and hence are
planar pencils of lines. The corresponding focal
hypersurfaces $M_i$ degenerate into n submanifolds of codimension 2. In the
examples discussed below focal submanifolds $M_i$ are
glued together to form an algebraic variety $V^{n-1} \subset P^{n+1}$ of
codimension 2, so that the lines of the
congruence (\ref{I2}) can be
characterized as $n$-secants of $V^{n-1}.$ It may happen that $V^{n-1}$
contains a linear subspace
of codimension 2. This is closely
related to the property for system (\ref{cons}) to possess {\it Riemann
invariants}.

\medskip
\noindent{\bf Definition.} {\it The Riemann invariant for the $i$th
characteristic velocity $\lambda ^i$ is  a function
$R^i(u)$ such that 
$$
R^i_t=\lambda ^iR^i_x
$$ 
by virtue of (\ref{cons}).}

There is a simple criterion for the existence of Riemann invariants in term
of the coefficients $c^i_{jk}$ defined by
(\ref{comm}). 

\begin{proposition}
The characteristic velocity $\lambda^i$ possesses a Riemann invariant if and
only if $c^i_{jk}=0$ for any $(j,k)\ne i.$
\label{Rim}
\end{proposition}
\noindent 

\begin{theorem}
If the characteristic velocity $\lambda ^i$ of a T-system (\ref{cons})
possesses a Riemann invariant, then the corresponding focal submanifold
$M_i$ is a linear subspace
of codimension 2. 
\label{RimSurf}
\end{theorem}
\noindent {\it Proof:}
Let us consider, for definiteness, the focal hypersurface $M_1$ with the
radius vector 
$$
{\bf r}_1=(\lambda ^1, u^1\lambda ^1-f^1,..., u^n\lambda ^1-f^n),
$$
corresponding to the characteristic velocity $\lambda ^1.$
We will need the following relations between the densities $u$ and the
fluxes $f$ of conservation laws of system (\ref{cons}):
\begin{equation}
L_i(f)=\lambda^iL_i(u)~{\rm for~ any~}i=1,...,n,
\label{LfLu}
\end{equation}
\begin{equation}
L_iL_j(u)=\frac{L_j(\lambda
^i)}{\lambda^j-\lambda^i}L_i(u)+\frac{L_i(\lambda
^j)}{\lambda^i-\lambda^j}L_j(u)
+\frac{\lambda^j-\lambda^k}{\lambda^i-\lambda^j}c^k_{ij}L_k(u),~i\ne j,
\label{2d}
\end{equation}
(see e.g. \cite{Sevennec}, \cite{Tsarev}). In particular, $f=f^s$ and
$u=u^s$ satisfy (\ref{LfLu}) and (\ref{2d})
for any $s=1,...,n.$ Introducing ${\bf l}=(1, \ u^1, ..., u^n)$ and applying
$L_1,...,L_n$ to the radius vector ${\bf r}_1$,
one readily obtains
$L_1({\bf r}_1)=L_1(\lambda ^1){\bf l}=0$ as $L_1(\lambda^1)=0$ by the
linear degeneracy.
Thus, the condition
$L_1({\bf r}_1)=0$ implies that $M_1$ is independent of $R^1$, where $R^1$
is the Riemann invariant corresponding
to $\lambda^1$. Since
$$
L_i({\bf r}_1)=
L_i(\lambda ^1){\bf l}+(\lambda ^1-\lambda ^i)L_i({\bf l}),
$$
the tangent space $TM_1$  is spanned by $n-1$ vectors
$L_i(\lambda ^1){\bf l}+(\lambda ^1-\lambda ^i)L_i({\bf l}),~(i\ne 1).$
 This tangent space belongs to the hyperplane $H_1$ spanned by
n vectors ${\bf l}, L_i({\bf l}),~(i\ne 1)$ which
depends on the variable $R^1$ only since $L_k(H_1)\in H_1$ for any $k\ne 1$.
The latter
follows from the relations
$$
\L_k({\bf l})\in H_1,~~~ L_k^2({\bf l})=A_kL_k({\bf l})\in H_1, ~~~
L_kL_jH_1\in H_1 ~~ (k, j \ne 1).
$$
Here $L_k^2({\bf l})=A_kL_k({\bf l})$ due to the linearity of rarefaction
curves, and $L_kL_jH_1\in H_1 $ by virtue of
(\ref{2d}) and the condition $c^1_{kj}=0$ (Proposition 1).
On the other hand, $M_1$ (and hence $TM_1$) does not depend on $R^1$ due to
the
linear degeneracy.  Consequently, the tangent space of  $M_1$ is the
intersection of any
two hyperplanes $H_1$ which
correspond to the
two different values of $R^1$. Therefore,
it is stationary and coincides with the focal submanifold $M_1$. Q.E.D.

\medskip

Finally, we recall the necessary information about {\it reciprocal
transformations} of systems
of conservation laws.
Let $B(u)dx+A(u)dt$ and $N(u)dx+M(u)dt$ be two  conservation laws of system
(\ref{cons}), understood as
the one-forms which are closed by virtue of (\ref{cons}).
 In the new independent variables $X,T$ defined by
\begin{equation}
dX=B(u)dx+A(u)dt,~~~dT=N(u)dx+M(u)dt,
\label{reciprocal}
\end{equation}
 system (\ref{cons}) takes the form
\begin{equation}
u^T_t=V^i_j(u)  u^j_X, ~~i=1,...,n,
\label{3hr}
\end{equation}
where $V=(Bv-AE)(ME-Nv)^{-1},~E=id.$
The new characteristic velocities $\Lambda ^k$ are
\begin{equation}
\Lambda ^k=\frac{\lambda ^kB-A}{M-\lambda ^kN}.
\label{newEigen}
\end{equation}
Transformations (\ref{reciprocal}) are called {\it reciprocal}.
Reciprocal transformations are known to preserve the linear degeneracy (see
\cite{Fer4}).
Moreover, if both integrals (\ref{reciprocal}) are linear combinations of
the canonical 
integrals $u^idx+f^idt$ defining  the
T-system (\ref{cons}),
\begin{equation}
\begin{array}{c}
dX=(\alpha _iu^i+\alpha) dx+(\alpha _if^i+\tilde{\alpha}) dt,\\
dT=(\beta _iu^i+\beta) dx+(\beta _if^i+\tilde{\beta} ) dt,
\end{array}
\label{Trec}
\end{equation}
(here $\alpha _i,\alpha,\tilde{\alpha},\beta _i,\beta,\tilde{\beta}$ are
arbitrary constants),
then the transformed system will be a T-system, too (\cite{Fer1}).
Futhermore, affine transformations
\begin{equation}
U^i=C^i_ju^j+D^i,~~C^i_j=const,~D^i=const,~{\rm det \ C^i_j \ne 0},
\label{affine}
\end{equation}
obviously transform T-systems to  T-systems.

\begin{theorem}\cite{Fer1}
The transformation group generated by reciprocal transformations
(\ref{Trec}) and affine transformations (\ref{affine})
is isomorphic to the group of projective
transformations of $P^{n+1}$.
\label{projective}
\end{theorem}

Thus, the classification of systems of conservation laws up to
transformations 
(\ref{Trec}) and (\ref{affine}) is equivalent to the classification of the
corresponding congruences
up to projective equivalence. Actually, this observation was the main reason
for introducing
the geometric correspondence discussed in section 2.

\section{Geometry of the  examples}

\noindent Let us first recall some of the well-known properties of the
Verones\'e variety $V^2\subset P^5$
realising  $P^5$ as the space of  $3\times 3$ symmetric matrices
$Z^{ij},~i,j=0,1,2.$ Verones\'e variety
$V^2$ is the variety of matrices of rank 1
$$
Z=\left(
\begin{array}{ccc}
X^0X^0 & X^0X^1 &X^0X^2 \\
X^1X^0 & X^1X^1 &X^1X^2 \\
X^2X^0 & X^2X^1 &X^2X^2
\end{array}
\right).
$$
It can be viewed as the canonical embedding
$F: P^2 \rightarrow V^2\subset P^5$ defined by
\begin{equation}
Z^{ij}=X^iX^j,~i=0,1,2,
\label{Ver}
\end{equation}
where $X^0:X^1:X^2$ are homogeneous coordinates in $P^2$.
Verones\'e variety coincides with the singular locus of the {\it cubic
symmetroid}
defined by the equation
$$
{\rm det} \ Z^{ij}=0,
$$
which  also is the bisecant variety of $V^2$  consisting of symmetric
matrices of rank two.
Under the embedding (\ref{Ver}) each line in $P^2$ is mapped onto a conic
on $V^2$, therefore, Verones\'e variety carries a 2-parameter family of
conics. 
The projective automorphism group of $V^2$ coincides with the natural action
of $PSL_3$ on $P^5$
\begin{equation}
Z \to g^TZg, ~~~ g\in PSL_3,
\label{hom}
\end{equation}
which obviously preserves $V^2$.

Below we discuss in some more detail the geometry of congruences associated
with the equations 
(\ref{ass}), (\ref{eq1}) -- (\ref{eq6}).

\subsection{ Geometry of the equations  with no Riemann invariants}

In this subsection we discuss equations (\ref{ass}), (\ref{eq1}) and
(\ref{eq2}).
Rewritten as systems of conservation laws, they  do not possess Riemann
invariants, so that
the corresponding focal varieties will be irreducible. We explicitly
demonstrate that they concide
with different 
non-singular projections of the Verones\'e variety.

\bigskip

\noindent {\bf Equation (\ref{ass}).}  The focal variety of the
corresponding congruence
(\ref{I24}) 
is defined by (\ref{cuspidal})
\begin{equation}
y^0=\lambda, ~~~ y^1=a\lambda-b, ~~~ y^2=b\lambda-c, ~~~ y^3=c\lambda-b^2+ac
\label{focso21}
\end{equation}
where  $\lambda $ satisfies the characteristic equation
\begin{equation}
\lambda ^3+a\lambda^2-2b\lambda+c=0.
\label{charso21}
\end{equation}
One can verify that the three focal surfaces  (\ref{focso21})
corresponding to the three different values of $\lambda $ are, in fact,
"glued" together to form the
algebraic variety
defined in this affine chart by a system of seven cubics
\begin{equation}
\begin{array}{l}
(y^0)^3+y^0y^1-y^2=0, ~~~ (y^2)^2+y^3(y^0)^2=0, ~~~
y^1y^2y^3+y^0(y^3)^2-(y^2)^3=0, \\
\ \\
y^2(y^0)^2+y^1y^2+y^0y^3=0, ~~~ (y^3)^2-y^1(y^2)^2+y^0y^2y^3+y^3(y^1)^2=0,
\\
\ \\
y^0y^1y^3-y^0(y^2)^2-y^2y^3=0, ~~~ y^0y^2+y^1(y^0)^2+(y^1)^2+y^3=0.
\end{array}
\label{yfocso21}
\end{equation}
Variety (\ref{yfocso21}) coincides with the projection of the Verones\'e
variety $V^2\subset P^5$
$$
y^0=\frac{Z^{02}}{Z^{22}},~~
y^1=\frac{Z^{12}-Z^{00}}{Z^{22}},~~
y^2=\frac{Z^{01}}{Z^{22}},~~
y^3=-\frac{Z^{11}}{Z^{22}}
$$
from the point 
\begin{equation}
\left(
\begin{tabular}{ccc}
$Z^{00}$ & 0 & 0\\
0 & 0 & $Z^{00}$\\
0 & $Z^{00}$ & 0
\end{tabular}
\right)
\label{ppoint1}
\end{equation}
into $P^4$. 
Notice that this point does not belong to the bisecant variety $S(V^2)$ and
hence the projection is non-singular (we would like to thank A.~Oblomkov
for clarifying the structure of  its ideal).
Here we list some of the main
properties of this projection which are, of course, well-known.

1. The manifold of trisecant lines of the focal variety (\ref{yfocso21}) is
three-dimensional. 

2. For a  point $p$ on the focal variety the set of trisecants passing
through $p$ 
forms a planar pencil with the vertex  $p.$

3. The intersection of the abovementioned planar pencil with the focal
variety 
consists of the point $p$ and a conic. Let us demonstrate this by a direct
calculation.
Since $(1,\lambda,\lambda^2)^T$ is the eigenvector of the system
(\ref{cons3}) corresponding to the eigenvalue
$\lambda$, the rarefaction curve passing through $p$ is given parametrically
by
$$
(a, \ b, \ c)+s(1, \ \lambda, \ \lambda^2)=(a+s, \ b+s\lambda, \
c+s\lambda^2),
$$
s being the parameter (recall that rarefaction curves are lines). The
corresponding  pencil of lines
$$
\begin{array}{l}
y^1=(a+s)y^0-(b+s\lambda ),\\
y^2=(b+s\lambda )y^0-(c+s\lambda^2 ),\\
y^3=(c+s\lambda^2)y^0-(b+s\lambda )^2+(a+s)(c+s\lambda^2 )
\end{array}
$$
belongs to the plane with parametric equations
$$
\begin{array}{l}
y^0=X,\\
y^1=aX-b+Y,\\
y^2=bX-c+\lambda Y,\\
y^3=cX-b^2+ac+\lambda^2 Y.
\end{array}
$$
It can be readily verified that the intersection of this plane with the
focal 
variety consists of the point $X=\lambda,~Y=0$
and the parabola $Y+X^2+(a+\lambda)X+\lambda^2+a\lambda-2b=0.$

\bigskip

\noindent {\bf Equation (\ref{eq1}).} Rewritten as a system of conservation
laws
\begin{equation}
a_t=b_x, ~~~ b_t=c_x, ~~~ c_t=((1+bc)/a)_x,
\label{so31}
\end{equation}
this equation is associated with the congruence
\begin{equation}
y^1=ay^0-b, ~~~ y^2=by^0-c, ~~~ y^3=cy^0-(1+bc)/a,
\label{congrso31}
\end{equation}
the focal variety of which is defined by (\ref{cuspidal})
\begin{equation}
y^0=\lambda, ~~~ y^1=a\lambda-b, ~~~ y^2=b\lambda-c, ~~~
y^3=c\lambda-(1+bc)/a
\label{focso211}
\end{equation}
where  $\lambda $ satisfies the characteristic equation
\begin{equation}
\lambda ^3-\frac{b}{a}\lambda^2-\frac{c}{a}\lambda+\frac{1+bc}{a^2}=0.
\label{charso211}
\end{equation}
One can verify that the three focal surfaces  (\ref{focso211})
corresponding to the three different values of $\lambda $ are
glued together to form the algebraic variety
defined in this affine chart by a system of cubics
\begin{equation}
\begin{array}{l}
1+y^0(y^1)^2+y^1y^2=0, ~~~ y^1(y^0)^2-y^3=0, ~~~
(y^0)^3+y^0y^2y^3+(y^3)^2=0, \\
\ \\
y^0+y^0y^1y^2+y^1y^3=0, ~~~ y^0y^3-y^2(y^0)^2+y^1(y^3)^2-y^3(y^2)^2=0, \\
\ \\
(y^0)^2+y^0y^1y^3+y^2y^3=0, ~~~ y^0y^1+(y^1)^2y^3-y^2-y^1(y^2)^2=0.
\end{array}
\label{yfocso211}
\end{equation}
This algebraic variety  is the projection
of the Verones\'e variety
$$
y^0=-\frac{Z^{02}}{Z^{12}},~~
y^1=-\frac{Z^{11}}{Z^{12}},~~
y^2=\frac{Z^{22}-Z^{01}}{Z^{12}},~~
y^3=-\frac{Z^{00}}{Z^{12}}
$$
from the point 
\begin{equation}
\left(
\begin{tabular}{ccc}
0 & $Z^{01}$  & 0\\
$Z^{01}$ & 0 & 0\\
0 & 0 &  $Z^{01}$ 
\end{tabular}
\right)
\label{ppoint2}
\end{equation}
into $P^4$.  
Notice that the two points (\ref{ppoint1}) and (\ref{ppoint2}) are
equivalent under the action of
the group (\ref{hom}) preserving the Verones\'e variety
(indeed, both matrices have the same Lorentzian signature). Hence,  both
projections and
the corresponding congruences of trisecants are projectively equivalent. To
be explicit, 
the projective transformation
\begin{equation}
y^0=-\frac{1}{Y^1}, \ y^1=\frac{Y^2}{Y^1}, \ y^2=\frac{Y^0}{Y^1}, \
y^3=-\frac{Y^3}{Y^1}
\label{transformation}
\end{equation}
identifies the systems of cubics (\ref{yfocso21}) and (\ref{yfocso211}).
Applying this transformation
to the congruence (\ref{I24}) and introducing the new parameters
$A=-1/c, \ B=b/c, \ C=a-b^2/c$, we readily rewrite (\ref{I24}) in the form
$$
Y^1=AY^0-B, ~~~ Y^2=BY^0-C, ~~~ Y^3=CY^0-(1+BC)/A
$$
which coincides with (\ref{congrso31}).
This gives geometric explanation of the transformation
between  equations (\ref{ass}) and (\ref{eq1}) mentioned in the
introduction. On the level of systems
of conservation laws (\ref{cons3}) and (\ref{so31}), this transformation is
a reciprocal equivalence.

\bigskip

\noindent {\bf Equation (\ref{eq2}).} Rewritten as a system of conservation
laws
\begin{equation}
a_t=b_x, ~~~ b_t=c_x, ~~~ c_t=((c^2+b^2-ac-1)/b)_x,
\label{so3}
\end{equation}
this equation is associated with the congruence
\begin{equation}
y^1=ay^0-b, ~~~ y^2=by^0-c, ~~~ y^3=cy^0-((c^2+b^2-ac-1)/b)
\label{congrso3}
\end{equation}
whose focal surfaces are glued together to form the algebraic variety which
is the projection 
of the Verones\'e variety
$$
y^0=-\frac{Z^{02}}{Z^{12}},~~
y^1=\frac{Z^{00}-Z^{22}}{Z^{12}},~~
y^2=\frac{Z^{01}}{Z^{12}},~~
y^3=\frac{Z^{11}-Z^{22}}{Z^{12}},
$$
from the point 
\begin{equation}
\left(
\begin{tabular}{ccc}
$Z^{00}$ & 0 & 0\\
0 & $Z^{00}$ & 0\\
0 & 0 &  $Z^{00}$ 
\end{tabular}
\right)
\label{point3}
\end{equation}
into $P^4$. 
Notice that this point is not equivalent (over the reals) to the points
(\ref{ppoint1}) and (\ref{ppoint2}) under the action of the group
(\ref{hom})  (indeed, the signature of (\ref{point3}) is Euclidean). Hence,
the 
congruence (\ref{congrso3}) is not projectively equivalent to any of the
congruences
(\ref{I24}) or (\ref{congrso31}). The corresponding systems of conservation
laws are 
not reciprocally related.

\bigskip

We point out that the Verones\'e variety $V^2\subset P^5$, being the
intersection of quadrics,
does not possess trisecant lines. Trisecants appear only after we project
$V^2$ into $P^4$.
Indeed, let $P_0$ be a point in $P^5$ not on the bisecant variety $S(V^2)$.
Viewed as a
$3\times 3$ symmetric matrix, $P_0$ defines a non-degenerate conic in $P^2$
\begin{equation}
\sum_{i,j=0}^{2}(P_0^{-1})^{ij}X^iX^j=0
\label{met}
\end{equation} 
where $X^0:X^1:X^2$ are homogeneous coordinates.
If a plane passes through $P_0$ and  cuts $V^2$ in three points, then
pre-images of these points under the embedding (\ref{Ver})
are pairwise conjugate
with respect to the conic (\ref{met}).
Conversely, $P_0$ lies in the  plane spanned by the images under (\ref{Ver})
of any three points  in $P^2$ which are pairwise conjugate
with respect to (\ref{met}).
Thus, there is a three-parameter family of  planes
passing through $P_0$ and cutting $V^2$ in three points.
Projecting this  family
from the point $P_0$ into $P^4$, we arrive at the congruence of lines in
$P^4$.
By a construction, its lines are trisecants of the projection
$\pi_{P_0}(V^2)$, which is thus the focal surface of
our congruence. To see that the developable surfaces of the congruence are
planar pencils of lines,
we consider a  line $L$ in $P^2$ defined by the equation
$L_0X^0+L_1X^1+L_2X^2=0$.
Under the embedding  (\ref{Ver}), this line corresponds to a
conic  on $V^2$  lying in the so called {\it conisecant plane} of $V^2$.
In matrix form equations of this plane are $LZ=0.$
The three-dimensional subspace $\Lambda$ spanned by $P_0$ and the conisecant
plane consists of 
all $Z$ such that
the vectors $LZ$ and $LP_0$ are collinear.
In addition to the conic in the conisecant plane, $\Lambda$ intersects $V^2$
in the point $P^L_0$
whose pre-image in $P^2$ under (\ref{Ver}) has homogeneous coordinates
$P_0L^T$. Consider now
the one-parameter family of planes in $P^5$ which belong to $\Lambda$ and
pass through the line joining
$P_0$ and $P^L _0$. Clearly, each of these planes  intersects $V^2$ in three
points, and 
the projection of this one-parameter family of planes into $P^4$
will be a planar pencil of lines. This gives developable surfaces of our
congruence.

\subsection{ Geometry of the equations with one Riemann invariant}

In this subsection we discuss equations (\ref{eq3}) and (\ref{eq4}).
Since both  equations possess only one Riemann invariant, the corresponding
focal 
varieties will be reducible, consisting of a cubic scroll and a plane
intersecting the cubic scroll
along its directrix.

\bigskip

\noindent {\bf Equation (\ref{eq3})} can be rewritten as a system of
conservation laws 
\begin{equation}
a_t=b_x, ~~ b_t=c_x, ~~ c_t=(bc/a)_x
\label{R1so21}
\end{equation}
the characteristic velocities of which are $\lambda ^1=b/a$ and
$ \lambda ^2,\lambda^3=\pm \sqrt{c/a}.$
The only Riemann invariant $R^1=c/a$  corresponds to $\lambda ^1$.
The focal surface corresponding to
$\lambda ^1$  is the  plane
\begin{equation}
y^1=y^3=0,
\label{R1plain}
\end{equation}
while the focal surfaces corresponding to $\lambda^2$ and $\lambda^3$
are glued together to form the cubic scroll defined by a system of quadrics
\begin{equation}
y^0y^1+y^2=0, ~~~ y^0y^2+y^3=0, ~~~ y^1y^3-(y^2)^2=0.
\label{yfocR1so21}
\end{equation}
The plane (\ref{R1plain}) intersects the cubic scroll along its directix
\begin{equation}
y^1=y^2=y^3=0.
\label{R1lineInt}
\end{equation}
The cubic scroll (\ref{yfocR1so21}) can be obtained by projecting the
Verones\'e variety
$$
y^0=\frac{Z^{02}}{Z^{12}},~~
y^1=\frac{Z^{11}}{Z^{12}},~~
y^2=-\frac{Z^{01}}{Z^{12}},~~
y^3=\frac{Z^{00}}{Z^{12}},
$$
from the point   
$$
\left(
\begin{tabular}{ccc}
 0  & 0 & 0\\
0 &  0  & 0\\
0 & 0 &  $Z^{22}$ 
\end{tabular}
\right).
$$
Notice that the center of this projection lies on the Verones\'e variety.
The directrix (\ref{R1lineInt}) is the image of the tangent plane
$Z^{00}=Z^{01}=Z^{11}=0$
to the Verones\'e variety in the centre of projection,
and the  plane (\ref{R1plain}) is the projection of the three-dimensional
linear
subspace in $P^5$ spanned by the tangent plane and the point
\begin{equation}
\left(
\begin{tabular}{ccc}
0   & $Z^{01}$ & 0\\
$Z^{01}$ &  0  & 0\\
0 & 0 &  0 
\end{tabular}
\right)
\label{point1}
\end{equation}
on the bisecant variety. Thus, the  focal variety of our congruence is
reducible and consists of
the plane (\ref{R1plain}) and the cubic scroll (\ref{yfocR1so21}).
Like in the case of systems without Riemann invariants,

1. the manifold of trisecants of the focal variety is three-dimensional, and

2. for a fixed point $p$ on the focal variety the set of trisecants passing
through $p$ 
forms a planar pencil with the vertex  $p$.
If $p$ belongs to the plane (\ref{R1plain}),
the corresponding planar pencil cuts the focal variety  in the point $p$ and
a conic. 
If $p$ belongs to the cubic scroll, it cuts the focal variety in the point
$p$ and a pair of
 lines. 

\bigskip

\noindent{\bf Equation (\ref{eq4})} can be rewritten as a system of
conservation laws 
\begin{equation}
a_t=b_x, ~~ b_t=c_x, ~~ c_t=((c^2+b^2-ac)/b)_x
\label{R1so3}
\end{equation}
the characteristic velocities of which are $\lambda ^1=c/b$ and
$\lambda ^2,\lambda^3=(c-a\pm \sqrt{4b^2+(c-a)^2})/2b.$
The only Riemann invariant $R^1=(c-a)/b$  corresponds
to $\lambda ^1$. The focal surface corresponding to
$\lambda^1$ is the  plane
\begin{equation}
y^1=y^3,~~y^2=0,
\label{R1plainso3}
\end{equation}
while the focal surfaces  corresponding to $\lambda^2$ and $\lambda^3$ are
glued together to form
the cubic scroll
defined by a system of quadrics
\begin{equation}
y^0y^3+y^2=0, ~~~ y^0y^2+y^1=0, ~~~ y^1y^3-(y^2)^2=0.
\label{yfocR1so3}
\end{equation}
The plane (\ref{R1plainso3})
intersects the cubic scroll (\ref{yfocR1so3}) along its directrix
\begin{equation}
y^1=y^2=y^3=0.
\label{R1lineIntso3}
\end{equation}
The cubic scroll (\ref{yfocR1so3}) can be obtained by projecting $V^2$
$$
y^0=\frac{Z^{02}}{Z^{12}},~~
y^1=\frac{Z^{00}}{Z^{12}},~~
y^2=-\frac{Z^{01}}{Z^{12}},~~
y^3=\frac{Z^{11}}{Z^{12}},
$$
from the point   
$$
\left(
\begin{tabular}{ccc}
 0  & 0 & 0\\
0 &  0  & 0\\
0 & 0 &  $Z^{22}$ 
\end{tabular}
\right)
$$
on $V^2$. 
The directrix (\ref{R1lineIntso3}) is the image of the
tangent plane $Z^{00}=Z^{01}=Z^{11}=0$ to $V^2$ in the centre of projection,
and the plane (\ref{R1plainso3}) is the image of the three-dimensional
linear subspace spanned by the tangent plane and the point
\begin{equation}
\left(
\begin{tabular}{ccc}
 $Z^{00}$  & 0 & 0\\
0 &  $Z^{00}$  & 0\\
0 & 0 &  0 
\end{tabular}
\right)
\label{point2}
\end{equation}
on the bisecant variety.

A coordinate-free construction of the congruences discussed above can be
described as follows.
Take a point $P_0 \in S(V^2)$ which is represented by a symmetric matrix of
rank two. Then there is
a nonzero vector $X\in P^2$ such that $P_0X=0$. Consider the tangent plane
to $V^2$ at the point
$F(X)$. The projection of $V^2$ into $P^4$ from the point $F(x)$ is a cubic
scroll. The projection
of the tangent plane is the directrix. The projection of the
three-dimensional space spanned by the tangent plane
and $P_0$ is the plane intersecting the cubic scroll along its directrix.

Although the last two examples look pretty similar,
they are not projectively equivalent, indeed, the points (\ref{point1}) and
(\ref{point2})
have different signatures.

\subsection{Geometry of the equation with two Riemann invariants.}

In this subsection we discuss equation (\ref{eq5}). Due to the existence of
two Riemann invariants,
the corresponding focal variety will be reducible consisting of two planes
and a two-dimensional
quadric.

\bigskip

\noindent{\bf Equation (\ref{eq5})} can be rewritten as a system of
conservation laws
\begin{equation}
a_t=b_x, ~~ b_t=c_x, ~~ c_t=((c^2-1)/b)_x
\label{R2}
\end{equation} 
with the characteristic velocities  $\lambda ^1=0$ and $\lambda
^2,\lambda^3=(c \mp 1)/b.$
The system has 2 Riemann invariants $(c\pm 1)/b$  corresponding to
$\lambda^2$ and  $\lambda^3$,
respectively.  The focal surfaces of the  associated congruence
corresponding 
to $\lambda^2$ and $\lambda^3$ are the planes
\begin{equation}
\begin{array}{l}
y^2=\mp 1,\\
y^0=\mp y^3,
\end{array}
\label{R2plains}
\end{equation}
while the third focal surface, corresponding to $\lambda ^1$, is the quadric
\begin{equation}
y^0=0,~~~y^1y^3-(y^2)^2+1=0.
\label{R2focsurf}
\end{equation}
The planes (\ref{R2plains}) intersect  the  quadric (\ref{R2focsurf}) along
the 
rectilinear generators
$$
y^0=0, ~ y^2=\mp 1, ~ y^3=0
$$
which belong to different families and meet at infinity.

One can describe this congruence in a coordinate-free form as follows.
Consider a quadric $Q$ in a hyperplane $\Lambda \subset P^4$.
Choose a point $p\in Q$ and draw two rectilinear generators $l_1, l_2$ of
$Q$ through $p$. 
Choose two  planes $\pi_1$ and $\pi_2$ which are not in $\Lambda$ such that
$l_i\subset \pi_i$ and
$\pi_1\cap \pi_2=p$. The union of $\pi_1,\ \pi_2$ and $Q$ is the focal
variety in question.
Its trisecants define a congruence in $P^4$.

\subsection{Geometry of the equation with three Riemann invariants.}

As follows from Theorem \ref{RimSurf}, the focal varieties of congruences
corresponding to  diagonalizable n-component T-systems are collections of n
linear subspaces of
codimension two in $P^{n+1}$. For $n=4$ we have 3 planes in $P^4$. To ensure
the nondegeneracy, 
we require that the points of their pairwise intersections are distinct.

\bigskip

\noindent{\bf Equation (\ref{eq6})} can be rewritten as a linear system of
conservation laws
\begin{equation}
a_t=b_x, ~~
b_t=c_x, ~~
c_t=b_x
\label{R3}
\end{equation} 
with the characteristic velocities $\lambda ^1=0,~\lambda
^2=1,~\lambda^3=-1.$
Being linear, this system has 3 Riemann invariants. The
focal surfaces of the associated conruence are the  planes
$$
y^1=y^3, ~y^0=0 ~~ {\rm for} ~~ \lambda^1=0,
$$
$$
y^3=-y^2, ~y^0=1 ~~ {\rm for} ~~ \lambda^2=1
$$
and
$$
y^3=y^2,~y^0=-1 ~~ {\rm for} ~~ \lambda^1=0,
$$
respectively.

\section{Linear congruences.}

A congruence (\ref{I2}) is called {\it  linear} (or general linear) if its
Pl\"ucker coordinates
$$
1, \ u^i, \ f^i, \ u^if^j-u^jf^i
$$
satisfy n linear equations of the form
\begin{equation}
\alpha +\alpha_iu^i+\beta_if^i+\alpha_{ij}(u^if^j-u^jf^i)=0
\label{linCongr}
\end{equation}
where $\alpha, \ \alpha_i, \ \beta_i, \ \alpha_{ij}$ are arbitrary constants
(notice that equations (\ref{linCongr}), being linear in $f$,
define $f^i$ as explicit functions of $u$). We emphasize that all
 examples discussed above belong to this class.

\begin{theorem}
\label{linearity}
Congruences corresponding to  three-component T-systems are linear.
\end{theorem}

\noindent The {\it Proof}  is technical and relegated to the Appendix.

\medskip 

\noindent
Let $q$ be a fixed point in $P^{n+1}$. For the  lines  of the congruence
(\ref{I2})
passing through $q$ we have
$f^i=u^iq^0-q^i$, which, upon the substitution into (\ref{linCongr}),
implies a linear system for $u$. In general, this system possesses a
unique solution, so that there exists a unique line of our congruence
passing through $q$ (such congruences are said to be of order one). The
focal variety $V$ (also called the jump locus) consists of those $q$ for
which the corresponding
linear system  is not uniquely solvable for $u$. One can show that
$V$ has codimension at least two, and in the case it equals two, the
developable surfaces
are planar pencils of lines. Moreover, the intersection of any of these
planes with
the focal variety $V$ consists of a point and a plane curve of degree $n-1$.

There are at least two different ways one could approach the classification
of three-component T-systems or, equivalently,
 the line congruences in $P^4$ whose developable surfaces are planar pencils
of lines. 
 The first way is to establish their linearity.
 In the parametrization (\ref{I2}) this means that the three-dimensional
surface with the
 radius-vector $(1, u, f, u\wedge f)$ representing our congruence in the
Grassmanian $G(1, 4)$
 lies in a linear subspace of codimension three. After the linearity is
established, the results of Castelnuovo
 \cite{Castelnuovo} (who demonstrated that the corresponding focal varieties
are
projections of the Veronese variety into $P^4$) complete the classification
and imply Theorem \ref{MainList}.
 
Another way makes use of the Theorem of Segre \cite{Segre2} saying that a
surface 
in projective space
carrying a two-parameter family of plane curves (not lines) is either a cone
or the 
surface of Veronese $V^2$ or its projection into $P^4$. Moreover, the
corresponding plane curves are conics.
This theorem is intimately related to our problem.
Indeed, let $M_1, M_2, M_3$ be three focal surfaces of our congruence in
$P^4$. 
Take a point $p\in M_1$ and consider the planar pencil of lines passing
through it. This plane
intersects $M_2$ and $M_3$ in the curves $\gamma_2$ and $\gamma_3$,
respectively. Varying $p$,
we conclude that both $M_2$ and $M_3$ contain two-parameter families of
plane curves and hence are projections
of the Verones\'e variety (the case of a cone can be easily ruled out).
Moreover, the curves $\gamma_2$ and $\gamma_3$ are conics. To show that both
$M_2$ and $M_3$ are actually parts
of one and the same Verones\'e variety, it is sufficient to demonstrate that
$\gamma_2$ amd $\gamma_3$ are
parts of one and the same plane conic. This can be done as follows:
intersect $\gamma_2$ and $\gamma_3$ by a line passing through $p$
and construct the tangent lines to $\gamma_2$ and $\gamma_3$ in the points
of intersection.
These tangent lines meet in a point $q$ lying in the same plane. Doing this
for all lines of the pencil with vertex $p$
we arrive at the curve $q$, which clearly must be a line (called the polar
of $p$) in case 
$\gamma_2$ and $\gamma_3$ are parts of one and the same conic.

Unfortunately, both proofs require differential identities which
do not immediately follow from the geometric data given. Thus, it proves
necessary 
to directly investigate the exterior differential system governing
three-component T-systems,
transforming it into the involutive form. Once it is done, the verification
of both properties mentioned above
reduces to a straightforward calculation.

In the case $n=4$ the geometry of focal varieties of general linear
congruences, known 
as the Palatini scrolls, was investigated in \cite{Palatini} (see also
\cite{Mezzetti} and
\cite{Ottaviani} for further properties of
the Palatini scrolls) (we thank F. Zak for providing these references).
T-systems of conservation laws corresponding to general linear congruences
will be discussed elsewhere.

\section{ Completely exeptional Monge-Amp\'ere type equations.}

Another important class of examples of T-systems is provided by completely
exceptional  Monge-Amp\'ere equations studied in \cite{Boi2}.
Equations of this type are defined as follows. Introduce the Hankel matrix
\begin{equation}
\left| 
\begin{array}{ccccc}
\frac{\partial ^{2m}u}{\partial x^{2m}}&\frac{\partial ^{2m}u}{\partial
x^{2m-1}\partial t}
&\frac{\partial ^{2m}u}{\partial x^{2m-2}\partial t^{2}} & ...
&\frac{\partial ^{2m}u}{\partial x^{m}\partial t^{m}}\\
\ \\
\frac{\partial ^{2m}u}{\partial x^{2m-1}\partial t}&
\frac{\partial ^{2m}u}{\partial x^{2m-2}\partial t^{2}}
&\frac{\partial ^{2m}u}{\partial x^{2m-3}\partial t^3}& ...
&\frac{\partial ^{2m}u}{\partial x^{m-1}\partial t^{m+1}}\\
\ \\
... &...&...&...&...\\
\ \\
\frac{\partial ^{2m}u}{\partial x^{m}\partial t^m}&
\frac{\partial ^{2m}u}{\partial x^{m-1}\partial t^{m+1}}
&\frac{\partial ^{2m}u}{\partial x^{m-2}\partial t^{m+2}}&...
&\frac{\partial ^{2m}u}{\partial t^{2m}}
\end{array}
\right|
\label{Hankel}
\end{equation}
and denote by $M_{J,K}(u)$ its minor of order $l$ whose rows and columns are
encoded in the multiindices
$J=(j_1,...,j_{l})$ and $K=(k_1,...,k_l)$, respectively.  PDE's in question
are
defined by  linear combinations of these minors,
\begin{equation}
\sum A^{J,K}M_{J,K}=0
\label{Boillat}
\end{equation}
where the summation is over all possible $l, J, K$, and $A^{J,K}$ are
arbitrary constants.
Any such equation can be rewritten as
$\frac{\partial ^{2m}u}{\partial t^{2m}}=f(\frac{\partial ^{2m}u}{\partial
x^{2m}},
\frac{\partial ^{2m}u}{\partial x^{2m-1}\partial t},...,
\frac{\partial ^{2m}u}{\partial x^{1}\partial t^{2m-1}})$, and
after the introduction of $a^1=\frac{\partial ^{2m}u}{\partial x^{2m}}, \
a^2=\frac{\partial ^{2m}u}{\partial x^{2m-1}\partial t}, \ ..., \
a^{2m}=\frac{\partial ^{2m}u}{\partial x^{1}\partial t^{2m-1}}$,
assumes the conservative form
\begin{equation}
a^1_t=a^2_x, ~~ a^2_t=a^3_x, ~~ ..., ~~ a^{2m}_t=f(a^1,a^2,...,a^{2m})_x.
\label{BoiSys}
\end{equation}
One can show that this is always a T-system (in fact, its linear degeneracy
was demonstrated in \cite{Boi2}),
and the corresponding congruence (\ref{I2}) has the following  properties:\\
-- its developable surfaces are planar pencils of lines,\\
-- its focal variety has codimension at least 2,\\
-- each developable surface intersects the focal variety in a point, which
is the vertex of the pencil,
and a plane curve of degree $n-1$.

To obtain systems of this type for odd $n$, one should consider  equations
(\ref{Boillat}) which are independent of
$\frac{\partial ^{2m}u}{\partial t^{2m}}$. Introducing $v=\frac{\partial
u}{\partial x}$ and rewriting the resulting
equation for 
$v$ (which is of the order $2m-1$) as a system of conservation laws, one
arrives at
the congruence (\ref{I2}) with the properties as formulated above. We are
planning to investigate the geometry of
these examples elsewhere. When $n\geq 4$ these congruences are not
necessarily linear.
In this case the focal varieties must be singular, as follows from
\cite{Mezzetti}.

\section{Appendix}

\subsection{The exterior representation of hydrodynamic type systems}

Investigating  nondiagonalizable  systems
\begin{equation}
u^i_t=v^i_j(u)  u^j_x, ~~i=1, 2, 3,
\label{3h}
\end{equation}
 it is convenient to use the following
exterior notation: let $l^i=(l^i_1(u),l^i_2(u),l^i_3(u))$ be left
eigenvectors of the matrix $v^i_j$ corresponding to the
eigenvalues $\lambda ^i$, i.e., $l^i_jv^j_k=\lambda ^i l^i_k.$ With the
eigenforms $\omega ^i=l^i_jdu^j$,
the system (\ref{3h}) is rewritten  in the following exterior form:
\begin{equation}
\omega ^i\wedge (dx+\lambda ^i dt)=0, ~~~~ i=1, 2, 3.
\label{exterior}
\end{equation}
Differentiation of $\omega ^i$ and $\lambda ^i$ gives the structure
equations 
\begin{equation}
d\omega ^i=-c^i_{jk}\omega ^j\wedge \omega^k, ~~~ (c^i_{jk}=-c^i_{kj}), ~~~
d\lambda ^i=\lambda ^i_j\omega ^j,
\label{structure}
\end{equation}   
containing all the necessary information about the system under study.
Notice that if $\omega^i(\xi_j)=\delta^i_j$ then the
 coefficients $c^i_{jk}$ are the same as that appearing in (\ref{comm}).

\subsection{The structure equations}

The three theorems formulated below show that the structure equations of
three-component T-systems
take surprisingly simple forms. We give the detailed proof
of the first theorem and only sketch the proofs of two others. Notice that
$\omega^i$ are defined up to a nonzero normalization $\omega^i \to p^i
\omega^i, \ p^i\ne 0$.

\begin{theorem}
The eigenforms of a  three-component  T-system with no Riemann invariants
can be  normalized so that the structure equations take the form
\begin{equation}
d\omega ^1=\omega ^2\wedge \omega ^3,~~~d\omega ^2=\epsilon \omega ^3\wedge
\omega ^1,~~~
d\omega ^3=\omega ^1\wedge \omega ^2,~{\rm where}~\epsilon=\pm 1.
\label{structureND}
\end{equation}
\end{theorem}
\noindent {\it Proof:} For a system with no Riemann invariants the forms
$\omega ^i$ can be normalized in such a way
that the structure equations take the form
\begin{equation}
\begin{array}{c}
d\omega ^1=a\omega ^1\wedge \omega ^2+b\omega ^1\wedge \omega ^3+\omega
^2\wedge \omega ^3,\\
d\omega ^2=p\omega ^2\wedge \omega ^1+q\omega
^2\wedge \omega ^3+\epsilon \omega ^3\wedge \omega ^1,\\
d\omega ^3=r\omega ^3\wedge \omega ^1+s\omega ^3\wedge \omega ^2+\omega
^1\wedge \omega ^2
\end{array}
\end{equation}
where $\epsilon =\pm 1$. Below we assume $\epsilon =1$, since the complex
normalization 
$$
\omega ^1\rightarrow i\omega ^1,\ \omega ^2\rightarrow \omega ^2,\ \omega
^3\rightarrow i\omega ^3
$$
reduces the case $\epsilon =-1$ to the case $\epsilon =1$ which allows to
treat both cases on  equal footing.
Since the systems under consideration are strictly hyperbolic, the
eigenforms $\omega ^i$ constitute a basis so that the
differential
 of a function $u$ can be decomposed as  $du=u_i\omega ^i$ where
$u_i=L_i(u)$.  
Differentiating the relations $du=u_i\omega ^i$, $df=\lambda ^iu_i\omega ^i$
(compare with (\ref{LfLu}))
and equating to zero  coefficients at
$\omega ^i\wedge \omega ^j$, one obtains

$$ \begin{array}{c} u_{12}=u_2 \frac{\lambda _1^2}{\lambda
^1-\lambda ^2}+u_1\frac{\lambda _2^1}{\lambda ^2-\lambda
^1}+u_1a+u_3\frac{\lambda ^3-\lambda ^2}{\lambda ^1-\lambda ^2}, \\
u_{21}=u_2 \frac{\lambda _1^2}{\lambda ^1-\lambda ^2}+u_1\frac{\lambda
_2^1}{\lambda ^2-\lambda ^1}+u_2p+u_3\frac{\lambda ^3-\lambda ^1}{\lambda
^1-\lambda ^2}, \\ u_{31}=u_3 \frac{\lambda _1^3}{\lambda ^1-\lambda
^3}+u_1\frac{\lambda _3^1}{\lambda ^3-\lambda ^1}+u_3r+u_2\frac{\lambda
^1-\lambda ^2}{\lambda ^1-\lambda ^3}, \\ u_{13}=u_3 \frac{\lambda
_1^3}{\lambda ^1-\lambda ^3}+u_1\frac{\lambda _3^1}{\lambda ^3-\lambda
^1}+u_1b+u_2\frac{\lambda ^3-\lambda ^2}{\lambda ^1-\lambda ^3}, \\
u_{32}=u_3 \frac{\lambda _2^3}{\lambda ^2-\lambda ^3}+u_2\frac{\lambda
_3^2}{\lambda ^3-\lambda ^2}+u_3s+u_1\frac{\lambda ^1-\lambda ^2}{\lambda
^2-\lambda ^3}, \\ u_{23}=u_3\frac{\lambda _2^3}{\lambda ^2-\lambda
^3}+u_2\frac{\lambda _3^2 }{\lambda ^3-\lambda ^2}+u_2q+u_1\frac{\lambda
^1-\lambda ^3}{\lambda ^2-\lambda ^3},  \end{array} $$
where $u_{ij}$ are defined as
$du_i=u_{ij}\omega ^j.$
  Let us introduce the vector $\stackrel{\rightarrow}{u}=(u^1,u^2,u^3)$. The
condition of the linearity of
  rarefaction curves in the
  coordinates $u^1, u^2, u^3$ can be written in the form
$$
\stackrel{\rightarrow}{ u_{11}}=A_1\stackrel{\rightarrow}{ u_1},
\ \ \stackrel{\rightarrow}{ u_{22}}=A_2\stackrel{\rightarrow}{ u_2},
\ \ \stackrel{\rightarrow}{ u_{33}}=A_3\stackrel{\rightarrow}{ u_3}
$$ 
where $A_1, A_2$ and $ A_3$ are certain proportionality factors.
Taking into account the known expressions for $\stackrel{\rightarrow}{
u_{ij}}$ when $i\ne j$, one obtains the following expressions
for $d\stackrel{\rightarrow}{ u_i}$:
$$
\begin{array}{l}
d\stackrel{\rightarrow}{u_1}=A_1{ \stackrel{\rightarrow}{u_1}}\,{
\omega^1}+\left ({ \stackrel{\rightarrow}{u_1}}\,\left (a+{\frac {{ \lambda
^1_2}}{{ \lambda ^2}-{ \lambda ^1}}}\right )+{\frac {{
\stackrel{\rightarrow}{u_2}}\,{ \lambda ^2_1}}{{ \lambda ^1}-{ \lambda
^2}}}+{\frac {{ \stackrel{\rightarrow}{u_3}}\,\left ({ \lambda ^3}-{ \lambda
 ^2}\right )}{{ \lambda ^1}-{ \lambda ^2}}}\right ){ \omega^2}\\ \quad
+\left
 ({ \stackrel{\rightarrow}{u_1}}\,\left (b+{\frac {{ \lambda ^1_3}}{{
\lambda ^3}- { \lambda ^1}}}\right )+{\frac {{
\stackrel{\rightarrow}{u_2}}\,\left ({ \lambda ^3}-{ \lambda ^2}\right )}{{
\lambda ^1}-{ \lambda ^3}}}+{\frac {{ \stackrel{\rightarrow}{u_3}}\,{
\lambda ^3_1}}{{ \lambda ^1}-{ \lambda ^3}}}\right ){ \omega^3} \end{array}
$$
\begin{equation} \label{dh} \begin{array}{l}
d\stackrel{\rightarrow}{u_2}=\left ({ \stackrel{\rightarrow}{u_2}}\,\left
(p+{\frac {{ \lambda ^2_1}}{{ \lambda ^1}-{ \lambda ^2}}}\right )+{\frac {{
 \stackrel{\rightarrow}{u_1}}\,{ \lambda ^1_2}}{{ \lambda ^2}-{ \lambda
 ^1}}}+{\frac {{ \stackrel{\rightarrow}{u_3}}\,\left ({ \lambda ^3}-{
 \lambda ^1} \right )}{{ \lambda ^1}-{ \lambda ^2}}}\right ){ \omega
^1}+A_2{
\stackrel{\rightarrow}{u_2}}\,{ \omega ^2}\\ \quad +\left ({
\stackrel{\rightarrow}{u_2}}\,\left (q+{\frac {{ \lambda ^2_3}}{{ \lambda
^3}-{ \lambda ^2}}}\right )+{\frac {{ \stackrel{\rightarrow}{u_1}}\,\left ({
\lambda ^1}-{ \lambda ^3}\right )}{{ \lambda ^2}-{ \lambda ^3}}}+{\frac {{
 \stackrel{\rightarrow}{u_3}}\,{ \lambda ^3_2}}{{ \lambda ^2}-{ \lambda
 ^3}}}\right ){ \omega ^3} \end{array}
 \end{equation}
 $$ \begin{array}{l}
d\stackrel{\rightarrow}{u_3}=\left ({  \stackrel{\rightarrow}{u_3}}\,\left
(r+{\frac {{  \lambda ^3_1}}{{  \lambda ^1}-{ \lambda ^3}}}\right )+{\frac
{{
  \stackrel{\rightarrow}{u_1}}\,{  \lambda ^1_3}}{{  \lambda ^3}-{ \lambda
  ^1}}}+{\frac {{  \stackrel{\rightarrow}{u_2}}\,\left ({  \lambda ^1}-{
  \lambda ^2} \right )}{{  \lambda ^1}-{  \lambda ^3}}}\right ){  \omega
^1}\\ \quad +\left ({ \stackrel{\rightarrow}{u_3}}\,\left (s+{\frac {{
\lambda ^3_2}}{{  \lambda ^2}-{  \lambda ^3}}} \right )+{\frac {{
 \stackrel{\rightarrow}{u_1}}\,\left ({  \lambda ^1}-{ \lambda ^2}\right
)}{{ \lambda ^2}-{  \lambda ^3}}}+{\frac {{  \stackrel{\rightarrow}{u_2}}\,{
\lambda ^2_3}}{{ \lambda ^3}-{  \lambda ^2}}}\right ){  \omega ^2}+A_3{
\stackrel{\rightarrow}{u_3}}\,{  \omega ^3}  \end{array} $$
where $\stackrel{\rightarrow}{u_i}=(u^1_i,u^2_i,u^3_i).$

Differentiating these equations and equating to zero  coefficients at
$\omega^i\wedge \omega ^j$, one obtains 9 equations
which are linear in $\stackrel{\rightarrow}{u_j}$. Since  $u^1, u^2, u^3$
are functionally independent, these equations
split with respect to $\stackrel{\rightarrow}{u_j}$, thus providing 27
equations for the derivatives of $\lambda ^k$,
the factors $A_1, A_2, A_3$ and
the coefficients of the structure equations $a,b,p,q,r,s$. The coefficients
at 
$\stackrel{\rightarrow}{u_3}\omega^1\wedge \omega ^2$  and
$\stackrel{\rightarrow}{u_2}\omega^3\wedge \omega ^1$ of the result of
differentiation of the first equation
(\ref{dh}) allow 
to find $A_1$ and $\lambda ^1_1$. Similarly, the coefficients at
$\stackrel{\rightarrow}{u_3}\omega^1\wedge \omega ^2$  and
$\stackrel{\rightarrow}{u_1}\omega^2\wedge \omega ^3$ of the result of
differentiation of the second equation
(\ref{dh}) give  
$A_2$ and $\lambda ^2_2$. Finally, the coefficients at
$\stackrel{\rightarrow}{u_1}\omega ^2\wedge \omega ^3$  and
$\stackrel{\rightarrow}{u_2}\omega ^3\wedge \omega ^1$ of the result of
differentiation of the third equation  (\ref{dh}) give
$A_3$ and $\lambda ^3_3$,

\begin{equation}
\label{L11}
\lambda ^1_1={\frac {\left (2\,p-2\,r\right )\left ({ \lambda ^1}-{ \lambda
^3}
\right )\left ({ \lambda ^2}-{ \lambda ^1}\right )}{{ \lambda ^2}-{ \lambda
^3
}}},
\end{equation}
\begin{equation}
\label{L22}
\lambda ^2_2={\frac {\left (2\,a-2\,s\right )\left ({ \lambda ^2}-{ \lambda
^3}
\right )\left ({ \lambda ^1}-{ \lambda ^2}\right )}{{ \lambda ^1}-{ \lambda
^3
}}},
\end{equation}
\begin{equation}
\label{L33}
\lambda ^3_3={\frac {\left (2\,q-2\,b\right )\left ({ \lambda ^2}-{ \lambda
^3}
\right )\left ({ \lambda ^1}-{ \lambda ^3}\right )}{{ \lambda ^1}-{ \lambda
^2
}}},
\end{equation}

\begin{equation}
\label{A}
A_1={\frac {\left (r-p\right )\left ({ \lambda ^2}-2\,{ \lambda ^1}+{
 \lambda ^3}\right )}{{ \lambda ^2}-{ \lambda ^3}}}+{\frac {2\,{ \lambda
^3_1}
\,\left ({ \lambda ^2}-{ \lambda ^1}\right )}{\left ({ \lambda ^1}-{
\lambda ^3}\right )\left ({ \lambda ^2}-{ \lambda ^3}\right )}}+{\frac {2\,{
 \lambda ^2_1}\,\left ({ \lambda ^3}-{ \lambda ^1}\right )}{\left ({ \lambda
^2
}-{ \lambda ^3}\right )\left ({ \lambda ^2}-{ \lambda ^1}\right )}},
\end{equation}
\begin{equation}
\label{B}
A_2={\frac {\left (s-a\right )\left ({ \lambda ^1}-2\,{ \lambda ^2}+{
 \lambda ^3}\right )}{{ \lambda ^1}-{ \lambda ^3}}}+{\frac {2\,{ \lambda
^3_2}
\,\left ({ \lambda ^1}-{ \lambda ^2}\right )}{\left ({ \lambda ^1}-{
\lambda ^3}\right )\left ({ \lambda ^2}-{ \lambda ^3}\right )}}+{\frac {2\,{
 \lambda ^1_2}\,\left ({ \lambda ^3}-{ \lambda ^2}\right )}{\left ({ \lambda
^1
}-{ \lambda ^3}\right )\left ({ \lambda ^1}-{ \lambda ^2}\right )}},
\end{equation}
\begin{equation}
\label{C}
A_3={\frac {\left (q-b\right )\left ({ \lambda ^1}-2\,{ \lambda ^3}+{
 \lambda ^2}\right )}{{ \lambda ^1}-{ \lambda ^2}}}+{\frac {2\,{ \lambda
^1_3}
\,\left ({ \lambda ^2}-{ \lambda ^3}\right )}{\left ({ \lambda ^1}-{
\lambda ^3}\right )\left ({ \lambda ^1}-{ \lambda ^2}\right )}}+{\frac {2\,{
 \lambda ^2_3}\,\left ({ \lambda ^3}-{ \lambda ^1}\right )}{\left ({ \lambda
^2
}-{ \lambda ^3}\right )\left ({ \lambda ^1}-{ \lambda ^2}\right )}}.
\end{equation}
For linearly degenerate systems the first three of these equations imply
\begin{equation}
r=p,\ s=a, \ q=b.
\label{rsq}
\end{equation}
For the eigenvalues of linearly degenerate systems and the first derivatives
thereof we have
\begin{equation}
\begin{array}{l}
d\lambda ^1=\lambda ^1_2 \omega ^2+\lambda ^1_3 \omega ^3\\
d\lambda ^2=\lambda ^2_1 \omega ^1+\lambda ^2_3 \omega ^3\\
d\lambda ^3=\lambda ^3_1 \omega ^1+\lambda ^3_2 \omega ^2\\
\end{array}
\label{dL}
\end{equation}
\begin{equation}
\begin{array}{l}
d\lambda ^1_2=(p\lambda ^1_2-\lambda ^1_3) \omega ^1+\lambda ^1_{22} \omega
^2+\lambda ^1_{23}\omega ^3\\
d\lambda ^1_3=(r\lambda ^1_3+\lambda ^1_2) \omega ^1+(\lambda
^1_{23}-q\lambda^1_2+s\lambda ^1_3) \omega ^2+
\lambda ^1_{33}\omega ^3\\
d\lambda ^2_1= \lambda ^2_{11}\omega ^1+(a\lambda ^2_1+\lambda ^2_3) \omega
^2+ 
(\lambda ^2_{31}+b\lambda ^2_1-r\lambda ^2_3)\omega ^3\\
d\lambda ^2_3=\lambda ^2_{31} \omega ^1+ (s\lambda ^2_3-\lambda ^2_1)\omega
^2+ \lambda ^2_{33}\omega ^3\\
d\lambda ^3_1= \lambda ^3_{11}\omega ^1+ \lambda ^3_{12}\omega ^2+(b\lambda
^3_1-\lambda ^3_2) \omega ^3\\
d\lambda ^3_2= (\lambda ^3_{12}-a\lambda ^3_1+p\lambda ^3_2)\omega ^1+
\lambda ^3_{22}\omega ^2
+(q\lambda ^3_2+\lambda ^3_1) \omega ^3
\end{array}
\label{dL1}
\end{equation}
Taking into account (\ref{rsq}) and substituting the expressions for $A_1,
A_2$ and $A_3$ back into (\ref{dh}), one obtains a closed
system for $\stackrel{\rightarrow}{u_1}$, $\stackrel{\rightarrow}{u_2}$,
$\stackrel{\rightarrow}{u_3}$. Differentiating this system once again and
collecting coefficients at
$\stackrel{\rightarrow}{u_i}\omega ^j\wedge \omega ^k$, one arrives at the
27 compatibility conditions involving
the second derivatives of $\lambda^i$
and the first derivatives of $a,b$ and $p$. These conditions allow to find
all second derivatives of $\lambda^i$.
Moreover, they imply $a=b=p=0.$

Indeed,  differentiating the first equation (\ref{dh}) and collecting
coefficients at
$\stackrel{\rightarrow}{u_1}\omega ^1\wedge \omega ^2$,
$\stackrel{\rightarrow}{u_2}\omega ^1\wedge \omega ^2$,
$\stackrel{\rightarrow}{u_1}\omega ^2\wedge \omega ^3$,
$\stackrel{\rightarrow}{u_2}\omega ^2\wedge \omega ^3$ and
$\stackrel{\rightarrow}{u_3}\omega ^3\wedge \omega ^1$,
we  calculate
$\lambda ^3_{12}$,
$\lambda ^2_{11}$,
$\lambda ^1_{23}$, 
$\lambda ^2_{31}$ and
$\lambda ^3_{11}$,
respectively.
(We do not write down these intermediate expression for nonzero $a,b,p$ due
to their complexity.)
Substituting these derivatives back into
the coefficients at $\stackrel{\rightarrow}{u_3}\omega ^2\wedge \omega ^3$
and $\stackrel{\rightarrow}{u_1}\omega ^3\wedge \omega ^1$, one  arrives at
the equations  
$$
\begin{array}{l}
pa-5b-a_1=0\\
bp+5a-b_1=0
\end{array}
$$ 
which determine  $a_1$
 and $b_1$. It is remarkable that these equations for the structure
coefficients do not include the eigenvalues and
 the derivatives thereof.

Similarly, the differentiation of the second equation (\ref{dh}) defines
$\lambda ^1_{22}$,
$\lambda ^3_{22}$,
by calculating   
the coefficients at
$\stackrel{\rightarrow}{u_1}\omega ^1\wedge \omega ^2$ and
$\stackrel{\rightarrow}{u_3}\omega ^2\wedge \omega ^3$, respectively, along
with the equations
$$
\begin{array}{l}
pa+5b-p_2=0\\ 
ab-5p-b_2=0\\
bp+3a-p_3=0\\
ab-2a_3+b_2-p=0
\end{array}
$$ 
which result from the coefficients at
$\stackrel{\rightarrow}{u_2}\omega ^1\wedge \omega ^2$,
$\stackrel{\rightarrow}{u_1}\omega ^3\wedge \omega ^1$,
$\stackrel{\rightarrow}{u_2}\omega ^3\wedge \omega ^1$ and
$\stackrel{\rightarrow}{u_2}\omega ^2\wedge \omega ^3$, respectively.

Finally, differentiation of the third equation (\ref{dh}) defines
$\lambda ^2_{33}$ and
$\lambda ^1_{33}$
from  
the coefficients at
$\stackrel{\rightarrow}{u_2}\omega ^2\wedge \omega ^3$ and
$\stackrel{\rightarrow}{u_1}\omega ^3\wedge \omega ^1$,  and, eventually,
the equations
$b=0$,  $p=0$ and $a=0$
from the coefficients at
$\stackrel{\rightarrow}{u_3}\omega ^1\wedge \omega ^2$,
$\stackrel{\rightarrow}{u_3}\omega ^2\wedge \omega ^3$ and
$\stackrel{\rightarrow}{u_3}\omega ^3\wedge \omega ^1$.

For $a=b=p=0$ not only the system (\ref{dh}), but also the system
(\ref{dL1}) is in involution. The second derivatives
$\lambda ^i_{jk}$ and $\lambda ^i_{jj}$ can be obtained by cyclic
permutations from 
\begin{equation}
\begin{array}{c}
\lambda ^3 _{12}=\lambda ^3_1\lambda ^3_2\left( \frac{1}{\lambda ^3-\lambda
^1}+ \frac{1}{\lambda ^3-\lambda ^2}\right)+
 \lambda ^3_2\lambda ^2_1\left( \frac{1}{\lambda ^1-\lambda ^2}+
\frac{1}{\lambda ^2-\lambda ^3}\right)
-\lambda ^3_1\lambda ^1_2\left( \frac{1}{\lambda ^3-\lambda ^1}+
\frac{1}{\lambda ^1-\lambda ^2}\right)+\\
+\lambda ^2_3\left( \frac{\lambda ^1-\lambda ^3}{\lambda ^1-\lambda
^2}\right)^2
-\lambda ^1_3\left( \frac{\lambda ^2-\lambda ^3}{\lambda ^2-\lambda
^1}\right)^2
\end{array}
\label{ijk}
\end{equation}  

\begin{equation}
\begin{array}{l}
\lambda ^2 _{11}=-2\lambda ^2_1\lambda ^3_1\left( \frac{1}{\lambda
^2-\lambda ^3}+ \frac{1}{\lambda ^3-\lambda ^1}\right)+
 2 \frac{(\lambda ^2_1)^2}{\lambda ^2-\lambda ^3}
+2\frac{(\lambda ^2-\lambda ^3)(\lambda ^2-\lambda ^1)}{\lambda ^3-\lambda
^1}
\end{array}
\label{ijj}
\end{equation}  
Q.E.D.

\medskip

The differential form $\omega ^i$ is proportional to the differential of a
function 
$\omega ^i=p^idR^i$
if and only if $d\omega ^i\wedge \omega ^i=0$ (this is a special case of the
Frobenius theorem). Recall that the function
$R^i$ is called the Riemann invariant of the system (\ref{3h}).

\begin{theorem}
The eigenforms of a  three-component T-system with one Riemann invariant
can be  normalized so that the structure equations take the form
\begin{equation}
d\omega ^1=\epsilon \omega ^2\wedge \omega ^3,~~~d\omega ^2=\omega ^3\wedge
\omega ^1,
~~~\omega ^3=dR^3,~{\rm where}~\epsilon=\pm 1.
\label{structureR1}
\end{equation}
\end{theorem}

\noindent {\it Proof:} Let the system (\ref{3h}) have only one Riemann
invariant $R^3$: $dR^3=\omega ^3.$
Then the structure equations (\ref{3h})
 can be normalized as follows:
\begin{equation}
d\omega ^1=a\omega ^1\wedge \omega ^2+b\omega ^1\wedge dR^3+c\omega ^2\wedge
dR^3,~~~
d\omega ^2=p\omega ^2\wedge \omega ^1+q\omega ^2\wedge dR ^3+dR^3\wedge
\omega ^1.
\label{structureR1a}
\end{equation}
As in the case with no Riemann invariants, we can find
$\vec {u}_{ij}=\vec {U}_{ij}(\lambda^l,\lambda^l_m,\vec {u}_n,a,b,c,p,q)$,
$i\ne j.$
Now, the linear degeneracy  and the compatibility conditions for
\begin{equation}
\begin{array}{l}
d \vec {u}_{1}=A_1\vec {u}_1\omega ^1+\vec {U}_{12}\omega^2+\vec
{U}_{13}\omega^3,~~~\\
d\vec {u}_{2}=\vec {U}_{21}\omega^1+A_2\vec {u}_2\omega ^2+\vec
{U}_{23}\omega^3,~~~\\
d\vec {u}_{3}=\vec {U}_{31}\omega^1+\vec {U}_{32}\omega ^2+A_3\vec
{u}_3\omega^3
\end{array}
\label{du}
\end{equation}
imply
\begin{equation}
\begin{array}{c}
a_1=p_2,~~~q_1=b_1,~~~q_2=b_2,~~~a_3-b_2=pc+aq,~~~b_1-p_3=a-pb,\\
c_1=0,~~~c_2=0,~~~c_3=2c(q-b).
\end{array}
\label{com1}
\end{equation}
The first five formulae (\ref{com1}) are equivalent to $d(p\ \omega ^1+a\
\omega ^2+b\ dR^3)=0$ and $d((q-b)\ dR^3)=0.$
Define the functions $\phi $ and $\psi$ by
$$
-\frac{d\phi}{\phi}=p\ \omega ^1+a\ \omega ^2+b\
dR^3,~~~\frac{d\psi}{\psi}=(q-b)dR^3.
$$
It is clear that $\psi$ depends on $R^3$ only. Renormalize the forms and
introduce the new Riemann invariant
$\widehat{R}^3$ as follows:
$$
\widehat{\omega}^1=\frac{\omega ^1}{\phi},~~~\widehat{\omega
}^2=\frac{\psi}{\phi}\omega^2,~~~d\widehat{R}^3=\psi(R^3)dR^3.
$$
Thus renormalized forms satisfy
$$
d\widehat{\omega} ^1=\frac{c}{\psi^2} \widehat{\omega} ^2\wedge
d\widehat{R}^3,~~~
d\widehat{\omega} ^2=d\widehat{R}^3\wedge \widehat{\omega }^1.
$$
The last three equations (\ref{com1}) give $d(c/\psi ^2)=0,$ so $c/\psi ^2$
is constant.
One can always choose $\psi$ to guarantee $c/\psi ^2=\pm 1.$
With the structure equation (\ref{structureR1}), the system (\ref{du}) is in
involution. 
Q.E.D. 

\medskip

\begin{theorem}
The eigenforms of a three-component T-system with two Riemann invariants
can be  normalized so that the structure equations take the form
\begin{equation}
d\omega ^1=\omega ^2\wedge \omega ^3,~~~\omega ^2=dR^2,~~~ \omega ^3=dR^3.
\label{structureR2}
\end{equation}
\end{theorem}
\noindent {\it Proof:} One can normalize  $\omega ^1$ so that
$$
d\omega ^1=a\ \omega ^1\wedge dR^2+b\ \omega ^1\wedge dR^3+dR^2\wedge dR^3.
$$
Now, the compatibility conditions for (\ref{du})    imply
$$
a_1=a_3=b_1=b_2=0,
$$
which means that $a$ is a function of $R^2$ and $b$ is a function of $R^3$
only.
Define 
$$
\begin{array}{c}
d\widehat{R}^2=exp\left(\int
a(R^2)dR^2\right)dR^2,~~~d\widehat{R}^3=exp\left(\int b(R^3)dR^3\right)dR^3,
\\
\ \\
\widehat{\omega}^1=exp\left(\int a(R^2)dR^2+\int b(R^3)dR^3\right)\omega ^1.
\end{array}
$$
The so renormalized forms and the so defined Riemann invariants satisfy
$d\widehat{\omega} ^1=d\widehat{R} ^2\wedge d\widehat{R} ^3.$
As before, with the structure equation (\ref{structureR2}), the
corresponding system (\ref{du}) is in involution.
Q.E.D. 

\medskip 

\noindent {\bf Remark.} Equations (\ref{eq1}) and (\ref{eq2}) have the
structure equations
(\ref{structureND}) with $\epsilon =-1$ and $\epsilon =1$, respectively.
Equations (\ref{eq3}) and (\ref{eq4}) have the structure equations
(\ref{structureR1})
with $\epsilon =-1$ and  $\epsilon =1$, respectively.  Equation (\ref{eq5})
has the structure equations (\ref{structureR2}).

\subsection{Proof of the Theorem \ref{linearity}}

In the parametrization (\ref{I2}) the linearity of the congruence means that
the three-dimensional submanifold with the
 radius-vector ${\bf q}=(1,{\bf u}, {\bf f}, {\bf u}\wedge {\bf f})$
representing the congruence in the Grassmanian
 $G(1,4)$ lies in a linear subspace of codimension three
 (here ${\bf u}=(u^1,u^2,u^3)$ and ${\bf f}=(f^1,f^2,f^3)$).
 The osculating space of this submanifold is spanned by the vectors ${\bf
q}_i=L_i({\bf q})$,
 ${\bf q}_{ii}=L_{i}^2({\bf q})$ and ${\bf q}_{ij}=L_jL_i({\bf q})$. The
conditions $\lambda ^i_i=0$,
 ${\bf u}_{ii}=A_i{\bf u}_i$ and
 ${\bf f}_i=\lambda ^i{\bf u}_i$ imply
  ${\bf q}_i=(0,{\bf u}_i,\lambda^i{\bf u}_i,{\bf u}_i\wedge ({\bf
f}-\lambda^i{\bf u}))$
 and
 ${\bf q}_{ii}=A_i{\bf q}_i \equiv 0\ {\rm mod}({\bf q}_1,{\bf q}_2,{\bf
q}_3).$ 
 Using the first equation of (\ref{dh}), one has
 $$
 \begin{array}{c}
{\bf q}_{12}=L_2({\bf q}_1)=\frac{\lambda ^2_1}{\lambda ^1-\lambda ^2}{\bf
q}_2+
 \lambda^2_1(0,0,{\bf u}_2,{\bf u}\wedge {\bf u}_2)+ \frac{\lambda
^1_2}{\lambda ^2-\lambda ^1}{\bf q}_1+
\lambda ^1_2(0,0,{\bf u}_1,{\bf u}\wedge {\bf u}_1)+ \\
\frac{(\lambda ^3-\lambda
^2)(\lambda^1-\lambda^3)}{\lambda^1-\lambda^2}({\bf q}_3+(0,0,{\bf u}_3,{\bf
u}\wedge {\bf u}_3))
 +(\lambda ^2-\lambda ^1)(0,0,0,{\bf u}_1\wedge {\bf u}_2),
\end{array}
$$
which, ${\rm mod}  ({\bf q}_1,{\bf q}_2,{\bf q}_3)$, equals
$$
\begin{array}{c}
\tilde {\bf q}_{12}= \lambda^2_1(0,0,{\bf u}_2,{\bf u}\wedge {\bf u}_2)+
\lambda ^1_2(0,0,{\bf u}_1,{\bf u}\wedge {\bf u}_1)+
 \frac{(\lambda ^3-\lambda
^2)(\lambda^1-\lambda^3)}{\lambda^1-\lambda^2}(0,0,{\bf u}_3,{\bf u}\wedge
{\bf u}_3)+\\
 +(\lambda ^2-\lambda ^1)(0,0,0,{\bf u}_1\wedge {\bf u}_2).
\end{array}
$$
 Define $\tilde {\bf q}_{23}$ and
 $\tilde {\bf q}_{31}$ in a similar way. Thus, the  osculating space is
spanned by six vectors
 ${\bf q}_1,{\bf q}_2,{\bf q}_3,\tilde{\bf q}_{12},
 \tilde{\bf q}_{23},\tilde{\bf q}_{31}$. Using the formula (\ref{ijk}) (and
the ones obtained from
 it by cyclic permutations) one shows by a direct computation that
$L_3(\tilde {\bf q}_{12}) \equiv 0\ \ {\rm mod}
 (\tilde{\bf q}_{12},
 \tilde{\bf q}_{23},\tilde{\bf q}_{31})$. This relation along with
 ${\bf q}_{ii}\equiv 0\ \ {\rm mod}  ({\bf q}_1,{\bf q}_2,{\bf q}_3)$
implies that the osculating space is stationary, so
the three-dimensional submanifold in question lies in six-dimensional linear
subspace.
In the cases with one or two Riemann invariants the proofs are essentially
the same. 
Q.E.D.

\section{Acknowledgements}
It is a great pleasure to thank F.L.~Zak and A.~Oblomkov for the very
usefull discussions and references. This research
was supported by the EPSRC grant No Gr/N30941.

\end{document}